# The Use of Fractional Blocks to Improve Mathematics for Second Grade Elementary School Students at South Bangka Indonesia

Suwarno[1], Raditya Bayu Rahadian [2]
[1] Education and Culture Office of South Bangka, Islands of Bangka Belitung, Indonesia
[2] Education Management, Graduate Schools, Yogyakarta State University, Indonesia

**Abstract:-** This study aims to improve mathematics learning results in the second grade of elementary school, especially the matter of fractions. The background of the research is the low learning outcomes achieved by students due to the learning carried out does not pay attention to the level of student development. This research was classroom action research. The research subjects were 23 students. The data collection instrument used the observation sheet and documentation. Quantitative data analysis techniques with descriptive statistics were used in this study to compare conditions before and after the action was carried out. The results showed an increase in mathematics learning outcomes in the subject of fractions. Student motivation shown in the implementation of learning in the second cycle reached 3.91 (high category) with N-Gain reaching 46.49%. The level of completeness has reached 100%. The average value obtained by students in the second cycle was 80.2 and N-Gain reached 34.16%. This increase was due to the learning carried out using fractional blocks. These results illustrate that the use of fractional blocks can improve mathematics learning outcomes, especially in the subject of fractions.

## I. INTRODUCTION

Low-level elementary students are still at the stage of concrete operational development, where students are able to think simple logic (Vandewaetere & Clarebout, 2014:p.24), but still have difficulty applying abstract concepts without the help of concrete props/media (Budiningsih, 2017:p.28). Meanwhile, mathematics subject matter for the second grade of the elementary school includes various abstract concepts, one of which is fractions (Purnomosidi, 2017).

The fraction matter discussed in second grade of the elementary school includes the concept of fractions $\frac{1}{2}$, $\frac{1}{3}$, and $\frac{1}{4}$. Most of the teachers have discussed this material with the help of pictures. The teacher theoretically explains the concept of fractions by giving examples of pictures of fractions on the board. Sometimes the teacher also uses other examples of pictures of fractions in providing explanations. Even so, the achievement of student learning outcomes is still low. From the results of tests conducted on 23 students from the second grade of SD 9 Simpang Rimba after discussing fractions, 14 students did not reach a minimum score of 70, with a very low average score is 65.0.

Reflections carried out by the teacher found several causes for the low student learning outcomes on fractions material in second grade of SD 9 Simpang Rimba. Judging from Piaget's theory of cognitive development, elementary-age students are at the stage of concrete operational thinking who already have logical thinking skills, but only through concrete objects, so that all learning components need to be adjusted according to these abilities (Winch & Gingell, 2008:p.58). The concrete here does not only mean that the object exists (can be seen, held, etc.), but also concrete explains its meaning. If what is meant is part of the whole, then the objects used must be separated and inserted according to their wishes in order to understand the concept. Therefore, the use of two-dimensional image media does not help students much in understanding the material.

Teacher-centered discussion of material also causes students not to be moved to try and find the meaning of the material being studied (Kenedi, Helsa, Ariani, Zainil, & Hendri, 2019:p.70). In accordance with the characteristics of the material, conceptual understanding cannot be transferred from teacher to student but must be through the students' own discovery and awareness with the help of the teacher. In addition, students are only faced with abstract problems that do not help their understanding (Purnomo, Widowati, & Ulfah, 2019:p.58), such as repeating the notion of fractions without giving students the opportunity to understand the material through their own way of thinking which is assisted by visual aids adequate.

Based on the background, it can be concluded that the low student learning outcomes are caused by the use of teaching aids and the presentation of learning activities that are not able to present the meaning of the concept of material for students. An alternative that can be used to improve student learning outcomes is to use props that are able to accommodate student activities through active practice and not as passive listeners.

One of the props that can be used to solve this problem is a broken block. This visual aid is in the form of three-dimensional wooden blocks that can be joined together according to the fractional material being studied. These three-dimensional props are expected to be able to cut the abstractness of the material so that it makes it easier for students to understand.

Based on the problems that have been described, the action hypothesis proposed in this research was:

1. The use of fractional blocks can improve student motivations to study mathematics, especially the subject of fractions in second grade of elementary school at Bangka Selatan, Indonesia.
2. The use of fractional blocks can improve the mathematics learning results in the second grade of elementary school, especially the subject of fractions at Bangka Selatan, Indonesia.

## II. LITERATURE REVIEW

*A. Mathematics for Second Grade of Elementary School in Indonesia*

Mathematics is a science that studies abstract structures and patterns of relationships that exist in them (Subarinah, 2006:p.1). Meanwhile, according to Prihandoko (2006:p.9) mathematics is concerned with structures, relationships, and abstract concepts developed according to logical rules. This means that learning mathematics is essentially learning concepts, conceptual structures, and looking for relationships between concepts and structures. Learning mathematics is essentially learning concepts, conceptual structures, and looking for relationships between concepts and their structures (Kenedi et al., 2019; Phonapichat, Wongwanich, & Sujiva, 2014), and should be connected to realistic (verbal) contexts, stay close to children, and be relevant to society (Boonen & Jolles, 2015:p.1). This understanding illustrates that mathematics can be viewed as a way/effort to gain knowledge/truth about something. This method/effort involves deductive reasoning, really a concept or statement that is stated as a logical truth from the previous one so that the concepts between concepts or statements in mathematics are consistent.

One of the mathematics materials in the second grade of the elementary school in Indonesia is fractions (Purnomosidi, 2017). A fraction in principle states several parts of the same number of parts. All the number of equal parts together form a unit (unit). Two kinds of circumstances that need emphasis are the overall concept as a unit and the same concept (Muhsetyo, 2007:p.4.5). The truth of a concept or statement is obtained as a logical consequence of the previous truth so that the relationship between concepts or statements in mathematics is consistent. This is evidenced by the background of the discussion regarding the emergence of fraction calculations.

The fractions material in second grade of the elementary school includes an introduction to the concept of fractions ½, $^1/_3$, and ¼. This concept includes identifying the correct form of fractions from several images similar to the proposed fraction, explaining the meaning of the fraction in question, and grouping fractions according to the desired shape (Purnomosidi, 2017). This material is abstract enough for $2^{nd}$ grade students of elementary school so that the implementation of learning must be able to accommodate the level of student development through activities that can be practiced through the help of appropriate teaching aids. This is intended so that the abstractness of the material can be more concrete. The provision of fraction material for the second grade of elementary school students is mainly to foster reasoning power in solving complex problems in a simple way. It is expected that students will be mentally trained to compete and not give up easily.

The obstacle that students often experience in studying fractions is related to their cognitive balance, which has long been established that the numbers are intact. According to Muhsetyo (2007:p.4.5-4.8), these difficulties include:
1. Students do not know the meaning of fractions.
2. Students have difficulty understanding fractions of a similar value.
3. Students have difficulty comparing and sorting fractions.
4. Students have difficulty adding and subtracting the unequal denominated fractions.

The difficulties that have been raised require the teacher to be able to provide views and invite students to have a new paradigm, that sometimes a calculation cannot be completed by only referring to integer calculations. The delivery should link students' real experiences, making it easier for students to grasp the abstractness of the material in a tangible form that they can feel. In the next development, if students have mastered the concept, calculations can be done that involve abstract reasoning students.

*B. Fractional Blocks as a Tools for Learning*

As material in the form of abstract concepts, the presentation of mathematics learning should be juxtaposed with the use of media/teaching aids that will make it easier for students to capture the material presented. Especially for elementary students who generally still have limitations in managing and understanding abstract concepts.

For low-grade elementary school students, learning mathematics can be assisted by the use of teaching aids that act as manipulative materials in order to understand mathematics learning material (Saleh, Prahmana, Isa, & Murni, 2018). Learning aids must be able to be used to present, present, present, or explain learning material to students, where the tools themselves are not part of the lesson given. These learning aids can be used by students to be directly involved in learning mathematics (Muhsetyo, 2007:p.2.3). From this understanding, it can be understood that if presenting learning material as well as practicing it will be more striking to students, than if it is only spoken.

The positive thing that is obtained from the use of media/teaching aids in learning is certainly influenced by the accuracy of selecting the media/teaching aids used with the characteristics of students and learning materials. The advantages that can be listened to from the use of media/teaching aids in learning according to Muhsetyo (2007:p.2.4) include:
1. Interesting and not boring students.

2. Easy to understand because it is assisted by visualization that clarifies the description.
3. Long-lasting memory of students.
4. Able to involve students' activeness in learning.

Especially for elementary students, especially those in low-grade classes, the teaching aids used must not only be able to actively involve students in learning, but must ensure the safety of these students in their use while minimizing the level of damage to the props against use, made of strong but harmless materials. for students, for example, wood, cardboard, and cloth (Muhsetyo, 2007:p.2.20).

The props used in this research are fractional blocks. This teaching aid is made of wood so that it is safe and durable against use because the user is a second grade student of elementary school who sometimes still can't be careful. These props can be done while playing. Students can hold, separate, combine, sort, and other practical activities in order to understand mathematics learning material about fractions. The use of these props is expected to concretize abstract mathematical concepts and increase student concentration and motivation if done while playing (Chizary & Farhangi, 2017; Spector, Merrill, Elen, & Bishop, 2014). Examples of fractional block props as shown in the following figure:

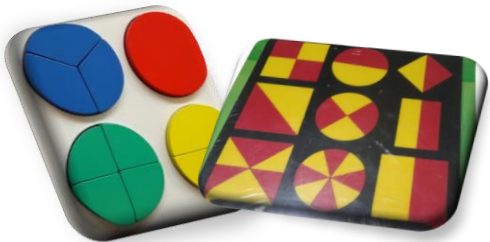

Fig 1. Examples of Fractional Blocks

*C. Student Learning Outcomes*

Learning outcomes are changes in behavior obtained by learning after experiencing learning activities (Budiningsih, 2012:p.14). While learning outcomes according to Sudjana are the abilities that students have after receiving their learning experiences (Sudjana, 2002:p.22). From the two definitions above, it can be concluded that learning outcomes are abilities or skills possessed by students after these students experience learning activities.

Gagne said there are five categories of learning outcomes: verbal information, intellectual skills, cognitive strategies, attitudes, and skills. Meanwhile, Bloom revealed three teaching goals which are one's abilities that must be achieved and are learning outcomes: cognitive, affective, and psychomotor (Sudjana, 2002:p.22). The learning outcomes achieved by students are influenced by two main factors, there are factors within students (internal) and from outside themselves (external). These factors can be detailed as follows (Thobroni, 2016:p.18):
1. Internal factors
- Physical (physiological) factors, both innate and acquired from the surrounding environment;
- Psychological factors, which consist of motivation, emotions, and others;
- Physical and psychological maturity factors.
2. External factors
- Social factors consisting of, family environment, school, community and group environment;
- Socio-cultural factors such as customs, science, technology, and the arts;
- Physical environmental factors: housing facilities, learning facilities, and climate;
- Religious spiritual environmental factors.

Learning outcomes achieved by students through the optimal teaching and learning process are indicated by the following characteristics (Sudjana, 2002:p.56):
1. Satisfaction and pride can foster intrinsic learning motivation in students. Students do not complain about low performance and they will strive harder to improve or at least maintain what has been achieved.
2. Increase his confidence and abilities, meaning that he knows his own abilities and believes that he has potential that is not inferior to others if he tries as he should.
3. The learning outcomes achieved are meaningful for him, such as being long-lasting, shaping behavior, useful for learning other aspects, willingness, and ability to develop creativity.
4. The learning outcomes obtained by students are comprehensive (comprehensive), which includes cognitive, knowledge or insight, affective (attitudes), and psychomotor, skills or behavior.
5. Students' ability to control or assess and control themselves, especially in assessing the results they achieve as well as assessing and controlling the process and learning efforts.

Learning outcomes are the terms used to express the level of success achieved by a person after going through the learning process. These three domains are the object of assessment of learning outcomes. Mathematics learning outcomes can be measured directly using learning outcomes tests. Learning outcome indicators help teachers and students measure progress and, if done correctly, helps engage and motivate students (Redfern, 2005:p.15). Based on the description that has been stated, it can be concluded that learning outcomes are an indicator of the success achieved by students in their learning efforts.

*D. Student Learning Motivations*

Motivation comes from the Latin word "movere" which means encouragement or directing. Motivation questions how to direct the power and potential of students so that they are willing to work together productively so as to achieve or achieve predetermined goals (Setiawan & Soeharto, 2020). The term motivation also comes from the word motive which can be interpreted as an effort to encourage someone to do something. Motive can be said as the driving force from within and within the subject to carry out certain activities in order to achieve a goal.

Motivation, when viewed in terms of human needs, consists of; 1) physiological (primary needs that must be satisfied), 2) security (inner, goods/objects), 3) social (feelings of love), 4) achievement (recognition from others), 5) self-actualization (Uno, Umar, & Panjaitan, 2014:p.61). Motivation is the power that drives someone to do an action. A student who has strong motivation will certainly study hard, and vice versa. Motivation is intrinsic (from within) and extrinsic (from outside). Teachers must be able to move student motivation from extrinsic to intrinsic (Setiawan & Soeharto, 2020) (Rachmawati & Daryanto, 2015).

The following principles must be met by learning to motivate students; 1) students must be actively involved in their own learning: 2) teachers must be flexible and make adjustments to teaching and learning according to student needs; 3) consideration must be given to the motivation and self-esteem of students as they become engaged in self-assessment (Middlewood, Richard, & Beere, 2005:p.171).

Based on the explanation that has been given, the researcher concludes that learning motivation is the energy that arises in students which becomes the impetus for achieving certain goals, namely achieving the desired learning outcomes. Furthermore, it can be said that to develop the motivation to learn properly, what must be done as a teacher are:
1. Designing or preparing interesting teaching materials.
2. Conditions for active learning.
3. Using fun learning methods and techniques.
4. Strive to meet the needs of students in learning (need to be respected, not feel pressured).
5. Convince students that they can make achievements.
6. Improve student work as soon as possible and convey the results as soon as possible.
7. Tells the story of the value of learning material achievement and relates it to the real-life of students every day.

## III. RESEARCH METHODOLOGY

### A. Research Design
This research was a classroom action research. Action research is a distinctive approach to inquiry that is directly relevant to classroom instruction and learning and provides the means for teachers to enhance their teaching and improve student learning (Stringer, 2014:p.1). This study uses the level 2 action research model in which the researcher does not conduct research to find problems and potentials, but directly tests actions that are believed to solve problems or can increase work effectiveness and efficiency (Sugiyono, 2015:p.53). The steps are as shown in the following figure:

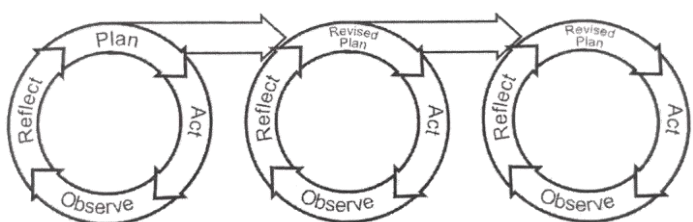

Fig 2. The Test of Action Hypothesis Process

Fig 2. shows that the action hypothesis testing is carried out using several cycles with the steps: planning, implementing, observing, and reflecting. Activities in testing the action hypothesis are carried out through several cycles with the following considerations (Sugiyono, 2015:p.55):
1. If the first cycle has been successful, then continue with the second cycle to test the consistency of the actions taken. If the test results in the first cycle are not different from the second cycle, then two cycles of testing can be sufficient.
2. If the test in the first cycle has not been successful, then after reflection is carried out, the planning is then revised, and testing is carried out in the second cycle. If the test in the second cycle has been successful, then the consistency is tested in the third cycle. If the results of the action test in the second and third cycles do not differ, then the action test can be terminated.

What is being experimented on or tried out is a plan of action or a hypothesis of action. Testing the action hypothesis using pre-experimental design (non-designs) in the form of One-Group Pretest-Posttest Design to compare conditions before and after being given treatment/action. With this form, the results of action can be known accurately (Sugiyono, 2015:p.137-138).

The treatment of this research is "the use of fractional blocks props". The experiment is "before-after", which is to compare student learning outcomes before the use of fractional blocks props and after using fractional block props. Success indicators are used as a reference in this study if the learning outcomes have met the following:
1. Student motivation is in the high category.
2. The average value of student learning outcomes reached 70.
3. Students' learning completeness is at least 80%.

### B. Population
This research involved students of second grade at elementary school 9 Simpang Rimba, South Bangka, Indonesia in the 2019/2020 academic year. The number of students who were the subject of this study amounted to 23 people, consisting of 11 female students and 12 male students.

### C. Instruments
The instruments used to observe changes in teacher work discipline are as follows. The score range refers to the five Likert scale: 5 = very good; 4 = good; 3 = medium; 2 = bad; 1 = very bad. In action research, the Likert scale is used to develop instruments to measure the attitudes, perceptions, and opinions of a person or group of people on the potential and problems of an object, action planning, and action results (Stringer, 2014:p.84).

Indicators of student activeness can be seen in several aspects, there are feelings of fun, attention, and interest in learning. The following are an explanation of the indicators for each aspect of student activity in learning (Sudjana, 2000:p.58):

#### 1. Feelings of fun
Fun feelings are shown by a cheerful face, smiling, and a friendly attitude in learning. The teacher notes the symptoms of feeling happy that the students show and gives a score of 1-5 according to the intensity of the occurrence.

*2. Attentions*

Attentions are shown by focusing on explanations, turning faces towards the teacher, not joking with friends, not making noise/talking alone. The teacher records the attention symptoms shown by students and gives a score of 1-5 according to the intensity of their occurrence.

*3. Interest in learning*

The students' interest was shown in learning such as taking notes, asking questions, answering questions, doing questions, and trying out props. The teacher records the students' symptoms of attraction and gives a score of 1-5 according to the intensity of their occurrence.

*D. Data Collection*

Data collection techniques used in each cycle of this research are observation and documentation. Observations were made during the implementation of learning by utilizing fractional blocks as props. When the learning process takes place, it will be observed how active students are in participating in learning activities. The documentation technique is used to obtain the achievement of learning outcomes before and after the action is taken.

*E. Data Analysis*

The data analysis technique used to test the hypothesis in each action cycle as follows:

1. To test the first action hypothesis, descriptive statistical analysis was used by means of a comparison of the average student learning motivations before and after the use of fractional block props. Student learning motivations are analyzed through a scale obtained from an assessment instrument using a Likert scale. To determine the level of student learning motivations using scores: Very High = 5, High = 4, Intermediate = 3, Low = 2, Very Low = 1 (Best & Kahn, 2006:p.331). The score is converted into a value with a five scale so that the assessment criteria are obtained as presented in the following table:

Table 1. Likert Scale Conversion

| Degree | Score Interval | |
|---|---|---|
| Very High (VH) | $X > \overline{X_i} + 1{,}80\ SBi$ | $X > 4{,}21$ |
| High (H) | $\overline{X_i} + 0{,}60\ SBi < X$ $\leq \overline{X_i} + 1{,}80\ SBi$ | $3{,}4 < X \leq 4{,}21$ |
| Intermediate (I) | $\overline{X_i} - 0{,}60\ SBi < X$ $\leq \overline{X_i} + 0{,}60\ SBi$ | $2{,}59 < X \leq 3{,}4$ |
| Low (L) | $\overline{X_i} - 1{,}80\ SBi < X$ $\leq \overline{X_i} - 0{,}60\ SBi$ | $1{,}79 < X \leq 2{,}59$ |
| Very Low (VL) | $X \leq \overline{X_i} - 1{,}80\ SBi$ | $X \leq 1{,}79$ |

From the assessment criteria shown in Table 1, it is obtained the standard level of student learning motivtions with the following details:

- Very High (VH) if the average score obtained is greater than 4.21.
- High (H) if the average score obtained is greater than 3.4 to 4.21.
- Intermediate (I) if the average score obtained is greater than 2.59 to 3.4.
- Low (L) if the average score obtained is greater than 1.79 to 2.59.
- Very Low (VL) if the average score obtained is less than 1.79.

2. To test the second action hypothesis, descriptive statistical analysis was used through a comparison of the average value of student learning outcomes. Quantitative analysis techniques with descriptive statistics can be explained as an analytical technique used to analyze data by describing or describing the collected data as it is without making generalized conclusions (Sugiyono, 2015:p.288).

## IV. RESULTS

Based on the description and data analysis that has been described, can be described the profile of the achievement of learning motivation of second grade students of SD 9 Simpang Rimba, South Bangka Regency before and after learning using fractional block props as presented in the following table:

Chart 1. The Development of Student's Learning Motivations

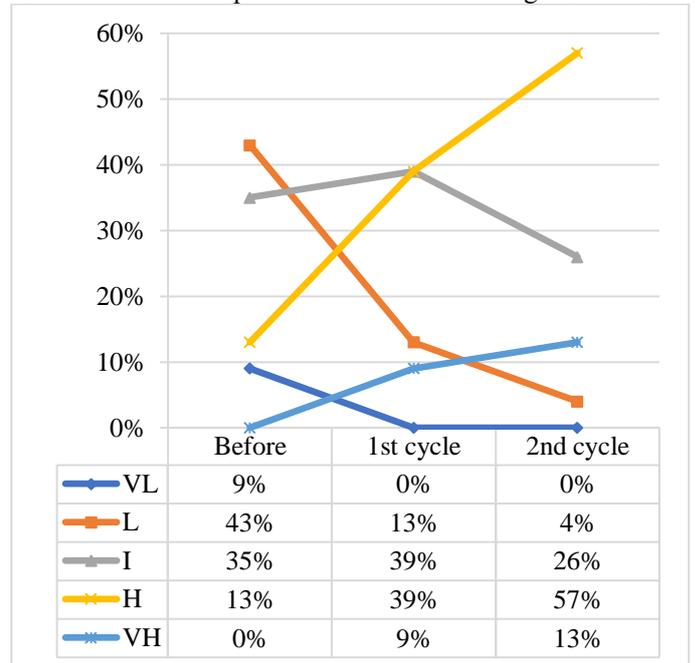

| | Before | 1st cycle | 2nd cycle |
|---|---|---|---|
| VL | 9% | 0% | 0% |
| L | 43% | 13% | 4% |
| I | 35% | 39% | 26% |
| H | 13% | 39% | 57% |
| VH | 0% | 9% | 13% |

Chart 1 shows the achievement of students' learning motivation if sorted according to level criteria. Before the implementation of the action, it was found that many students had low motivation. As actions are taken, students slowly begin to show active learning motivation. This was consistently shown until second cycle. The use of fractional blocks is able to motivate students to learn about fraction. The achievement of student learning motivations in detail are presented in the following table:

Table 2. Summary of Student Motivations Score's

| Number of cycles | Aspect of Student Learning Motivations | | | |
|---|---|---|---|---|
| | Feeling fun | Interest | Attention | Mean |
| 1st | 3.47 | 3.21 | 3.6 | 3.43 |
| 2nd | 3.95 | 3.82 | 3.95 | 3.91 |
| Mean | 3.71 | 3.52 | 3.78 | 3.67 |
| Before the action | 2.52 | 2.34 | 2.69 | 2.52 |

Table 2 provides information that the mean achievement of student learning motivations before the action is 2.52 (low). When observed after the use of fractional blocks the mean score in the first cycle increased to 3.43 (high), increased again in the second cycle to 3.91 (high). This means that student learning motivations have an increasing trend. As shown in the following chart:

Chart 2. Increase of Student Learning Motivations

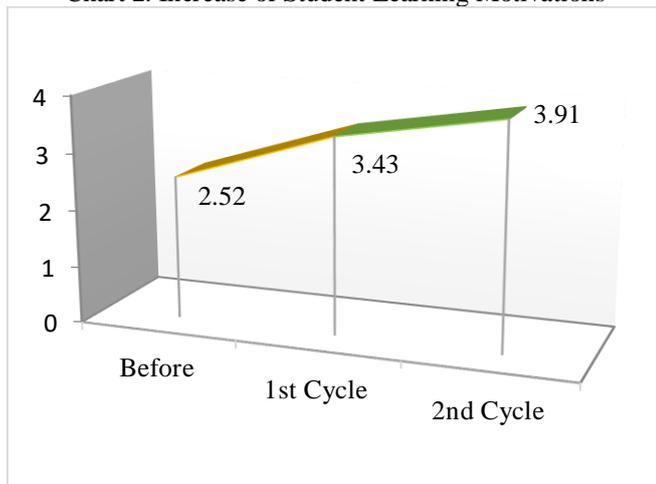

Chart 2 shows a significant increase in the achievement of student learning motivation between before and after taking action in the form of utilizing fraction blocks in learning mathematics in fractions material. Based on this graph, it can be interpreted that before the action is carried out, the majority of students have a low level of learning activeness. On an ongoing basis, in the application of each action cycle, student learning outcomes have increased steadily. This shows that the consistency of student learning activeness is maintained by the use of fraction blocks in learning mathematics especially the matter of fractions.

It cannot be denied that increasing student learning motivations directly affects students' interest in carrying out various learning activities in order to understand the learning material that has been prepared. This is what causes an increase in student learning outcomes as evidenced by the increase in the average value of student learning outcomes before and after the action. Details are presented in the following table:

Chart 3. Students Achievement Scores

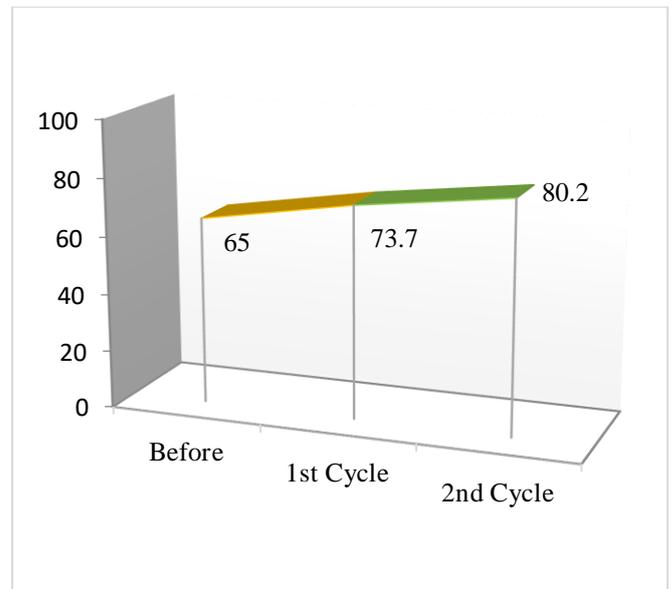

Chart 3 shows the consistency of improving student learning outcomes before and after taking action. This shows that the use of fraction blocks as a teaching aid in mathematics learning fractions material can consistently improve student learning outcomes. Before the action was taken, the student learning outcomes score only reached 65, still below the minimum criteria set at 70. As for the action in the first cycle, the average score of learning outcomes increased to 73.7, and in the second cycle continued to increase to 80.2. The increase in the average student learning outcomes shows an increase in student learning completeness as shown in the following graph:

Chart 4. Student Learning Completeness

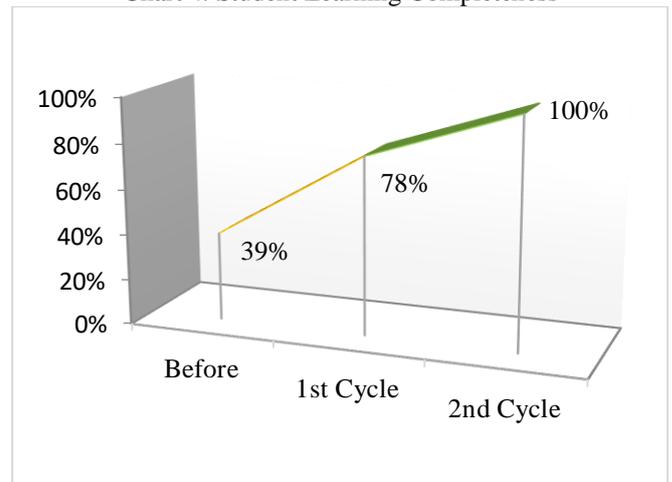

Related to the completeness of learning outcomes, Chart 4 informs an increase in the number of students who meet the minimum criteria set at 70 which is determined as a requirement for learning completeness. Before the action was taken, the students' learning completeness was 39%. In the first cycle, it increased to 78%, and in the second cycle, it increased again and even reached 100%.

## V. DISCUSSION

Daily activities cannot be separated from the role of mathematics. Starting from simple things when calculating household needs, as well as more complex needs related to the development of science. Learning mathematics is essentially learning concepts, conceptual structures, and looking for relationships between concepts and their structures (Kenedi et al., 2019; Phonapichat et al., 2014). Mathematics can be seen as a way to get the truth about something. This method involves deductive reasoning, in which the truth of a concept or statement is obtained as a logical result of the previous truth so that the relationship between concepts or statements in mathematics is consistent.

For elementary students, understanding mathematics with deductive reasoning is certainly very difficult, especially considering the inadequate ability to translate abstract concepts. So that learning and understanding of concepts can be initiated inductively through the experience of real events or intuition. Therefore mathematics learning in elementary school should emphasize the process of providing learning experiences to students through a series of planned activities so that the competence of the mathematics material being studied is achieved (Muhsetyo, 2007:p.1.26).

One of the matters in mathematics for elementary school is fractions. The provision of this matter is based on the complexity of the needs of human life which sometimes cannot be met by calculating using integers. Things that cannot be fulfilled by the existence of integers include stating some of the same parts of the whole, stating the number of several objects from a number of objects, calculating the percentage of parts, and much more (Muhsetyo, 2007:p.4.4). For elementary school students, giving material about fractions is mainly to foster reasoning power in solving complex problems in a simple way. It is expected that students will be mentally trained to compete and not give up easily.

Even so, in studying fractions of second grade of elementary school, students still often encounter difficulties regarding their cognitive balance which has long been established that the numbers are intact. This difficulty requires the teacher to be able to provide views and invite students to have a new paradigm, that sometimes a calculation cannot be completed by only referring to integer calculations. The delivery should link students' real experiences, making it easier for students to grasp the abstractness of the material in a tangible form that they can feel. In the next development, if students have mastered the concept, calculations can be done that involve abstract reasoning students.

For low-grade elementary school students, learning mathematics can be assisted by the use of teaching aids that act as manipulative materials in order to understand mathematics learning material (Saleh et al., 2018). This manipulative material serves to simplify difficult concepts, present relatively abstract materials to be more concrete, explain definitions or concepts more concretely, explain certain properties related to calculating (operations) and the properties of geometric shapes, and show the facts (Muhsetyo, 2007:p.2.20).

For elementary students, especially those in low-grade, the props used must not only be able to actively involve students in learning but also must ensure the safety of these students in their use while minimizing the level of damage to the props against use, made from strong materials but not hazardous to students, for example, wood, cardboard, and cloth (Muhsetyo, 2007:p.2.20).

The fractional block props used in this study are made of wood so that they are safe and durable against use because the users are second grade students of elementary school who sometimes still cannot be careful. These props can be done while playing. Students can hold, separate, combine, sort, and other practical activities in order to understand mathematics learning material about fractions. The use of this teaching aid is expected to be able to concretize abstract mathematical concepts and increase student concentration if done while playing (Chizary & Farhangi, 2017:p.237).

Implementation of learning using fractional block props follows the following procedure: 1) the teacher prepares fractional block props; 2) the teacher explains the concepts of fractions ½, $\frac{1}{3}$, and ¼ as well as demonstrating them using the fractional blocks; 3) students are divided into small groups; 4) each group is given a package of fractional blocks according to the concept of fractions ½, $\frac{1}{3}$, and ¼; 5) the teacher distributes student worksheets; 6) each student in the group practices the fraction blocks and fills in their worksheets; 7) the teacher checks the results of the student worksheets; 8) the teacher again explains the concepts of fractions ½, $\frac{1}{3}$, and ¼ as well as demonstrating them using the fraction blocks to strengthen students' understanding.

Students easily understand the concept of the material presented by practicing it using fractional blocks. In addition, students feel interested in learning because what is in front of them is a real object that can help them learn the material and not just an oral explanation. The division of small groups allows the intensity of the experiment done by each student to be more than if it was done classically. This further adds to the student's understanding of the fraction material because he can repeat it again and again until he really understands. Another interesting thing is that learning is not boring for students, as evidenced by the high learning motivation of students in trying, asking questions, and helping to explain to their friends.

The learning condition after this action is the opposite of what happened before the use of fractional blocks in learning. In the initial conditions, students are less interested in learning mathematics, which causes them not to participate in learning wholeheartedly. Low motivation so that learning outcomes are also low. This is indicated by the student's motivation score which reaches 2.52 which is in the low category. The average student learning outcomes reached 65.0 out of a scale of 100. Students' completeness only reached 39%. The low mathematics achievement of students is caused by learning that does not pay attention to the characteristics of the material and students. The abstract concept of mathematics is still difficult for students to understand if learning is not assisted by the use of props.

After the action, there was a consistent increase in student motivation. At the end of the second cycle, the students' motivation score reached 3.91 or was included in the high category. Student learning outcomes reached 80.2 and learning completeness reached 100%. This achievement occurred because the learning atmosphere was in favor of the students, caring about the students' difficulties in understanding the abstractness of the material using fraction teaching aids. Increased learning activities will also encourage students' willingness to study the material correctly, so as to improve their learning outcomes.

To determine the achievement of the effectiveness of learning using the use of fractional block props to increase student activity and learning outcomes, a comparison before and after the action (N-Gain) was used. N-Gain for student learning motivations reached 46.49% and N-Gain for student learning outcomes reached 34.16%. So, based on the discussion of the research results, it can be concluded that: 1) the use of fractional blocks can improve student motivations in mathematics, and 2) the use of fractional blocks can improve the mathematics results, in the second grade of elementary school especially fractions at SD 9 Simpang Rimba, South Bangka, Indonesia are proven empirical or acceptable.

## VI. CONCLUSIONS AND SUGGESTIONS

The use of fractional blocks can improve student motivations and learning results in mathematics especially the subject of fractions in the second grade of elementary school The steps for utilizing the fractional blocks in instructional follows the following procedure: 1) the teacher prepares fractional block props; 2) the teacher explains the concepts of fractions ½, $1/3$, and ¼ as well as demonstrating them using the fractional blocks; 3) students are divided into small groups; 4) each group is given a package of fractional blocks according to the concept of fractions ½, $1/3$, and ¼; 5) the teacher distributes student worksheets; 6) each student in the group practices the fraction blocks and fills in their worksheets; 7) the teacher checks the results of the student worksheets; 8) the teacher again explains the concepts of fractions ½, $1/3$, and ¼ as well as demonstrating them using the fraction blocks to strengthen students' understanding.

Fractional blocks can be used as learning tools to improve student motivations and learning outcomes in mathematics especially fractions. Even so, further research is needed to broadly enrich other practical steps in using the fractional blocks widely.

## ACKNOWLEDGMENT

The authors are thankful to Education and Culture Office of South Bangka, Islands of Bangka Belitung, Indonesia, for supporting this research work.